\documentclass[12pt]{amsart}
\usepackage{amsmath,amssymb}
\newtheorem{theorem}{Theorem}
\newtheorem{proposition}[theorem]{Proposition}
\newtheorem{lemma}[theorem]{Lemma}
\newtheorem{corollary}[theorem]{Corollary}

\def\eps{\varepsilon}
\def\ds{\displaystyle}

\title[Estimates of the Bergman distance on planar domains]
{Estimates of the Bergman distance on
Dini-smooth bounded planar domains}

\author{Nikolai Nikolov}

\author{Maria Trybu\l{}a}

\address{Institute of Mathematics and Informatics\\Bulgarian Academy
of Sciences\\ Acad. G. Bonchev 8, 1113 Sofia, Bulgaria\newline
\indent Faculty of Information Sciences\\
State University of Library Studies and Information Technologies\\
Shipchenski prohod 69A, 1574 Sofia, Bulgaria}\email{nik@math.bas.bg}
\address{Institute of Mathematics, Faculty of Mathematics and Computer Science,\\
Jagiellonian University, \L ojasiewicza 6, 30-348 Krak\'ow,
Poland}\email{maria.trybula@im.uj.edu.pl}


\thanks{M. Trybu\l{}a is partially supported by the International PhD programme
``Geometry and Topology in Physical Models'' of the Foundation for
Polish Science, by the Polish National Science Center -- grant PRO-2013/11/N/ST1/03609,
and by the Bulgarian National Science Found under contract DFNI-I 02/14.
The initial version of this paper was prepared during her visit to the
Institute of Mathematics and Informatics, Bulgarian Academy of
Science, October 2013 -- April 2014 and July 2014.}

\begin{document}

\begin{abstract} Precise estimates for the Bergman distances
of Dini-smooth bounded planar domains are given. These estimates
imply that on such domains the Bergman distance almost coincides
with the Carath\'eodory and Kobayashi distances.
\end{abstract}

\maketitle

\section{Results}

In \cite[Proposition 8]{Nikolov2}, the first named author found
optimal estimates for Carath\'{e}odory and Kobayashi distances,
$c_D$ and $k_D,$ on Dini-smooth bounded planar domains
$D$ in terms of $d_D=\mbox{dist}(\cdot,\partial D).$
In this paper we shall prove similar estimates for the Bergman
distance $b_D.$ For convenience of the reader, the definitions of
these three distances, as well as of Dini-smoothness, are given in
the next section.

\begin{proposition}
Let D be a Dini-smooth bounded planar domain. Then there exists a
constant $c>1$ such that
$$
\sqrt2\log\left(1+\frac{|z-w|}{c\sqrt{d_{D}(z)d_{D}(w)}}\right)\le
b_{D}(z,w)$$
$$\le\sqrt2\log\left(1+\frac{c|z-w|}{\sqrt{d_{D}(z)d_{D}(w)}}\right),
\quad z,\,w\in D.$$
\end{proposition}

By \cite[Proposition 8]{Nikolov2}, the same result holds for
$\sqrt2c_D$ and $\sqrt2k_D$ instead of $b_D.$ So, we have the
following

\begin{corollary} If $D$ is a Dini-smooth bounded planar domain,
then the differences $b_{D}-\sqrt{2}c_{D}$ and $b_{D}-\sqrt{2}k_{D}$
are bounded.
\end{corollary}

Note that Proposition 1 is equivalent to
\smallskip

\noindent{\bf Proposition 1$^\prime$.} {\it Let D be a Dini-smooth
bounded planar domain. There exists a constant $c>1$ such that:
\smallskip

\noindent$\bullet$ if $|z-w|^2>d_{D}(z)d_{D}(w),$ then
$$\log\frac{|z-w|^{2}}{d_{D}(z)d_{D}(w)}-c<\sqrt{2}
b_{D}(z,w)<\log\frac{|z-w|^{2}}{d_{D}(z)d_{D}(w)}+c;$$

\noindent$\bullet$ if $|z-w|^2\le d_{D}(z)d_{D}(w),$ then
$$\frac{|z-w|}{c\sqrt{d_{D}(z)d_{D}(w)}}\le
b_{D}(z,w)\le\frac{c|z-w|}{\sqrt{d_{D}(z)d_{D}(w)}}.$$}

\noindent{\bf Remark.} (a) The Dini-smoothness is essential as
an example of a $\mathcal{C}^{1}$-smooth bounded simply connected
planar domain shows (see \cite[Example 2]{Nikolov1}).
\smallskip

\noindent (b) One of the missing properties of $b_{D}$ in comparison
with $c_{D}$ and $l_{D}$ is monotonicity under inclusion of
(planar) domains. However, the invariants $M_{D}$ and $K_{D}$ share
this property which allows us to modify the approach from \cite{Nikolov2}.
\smallskip

\noindent (c) Results in $\mathbb{C}^{n}$ in the spirit of
Proposition 1 and Corollary 3 can be found in \cite{Bonk} and
\cite{Nikolov3}, respectively, where the strictly pseudoconvex
domains are treated. Note also that the Levi pseudoconvex corank one
domains are considered in \cite{Herbort}. As can be expected, our
estimates are more precise than those in \cite{Bonk} and \cite{Herbort},
when the two points, $z$ and $w,$ are close to each other.
\smallskip

\noindent (d) It follows by the second statement of Proposition
1$^\prime$ that
$$\frac{1}{cd_D(u)}\le\liminf_{\substack{z,w\to u\\z\neq
w}}\frac{b_D(z,w)}{|z-w|},\quad u\in D.$$ This inequality agrees
with the fact that (cf. \cite[Lemma 4.3.3 (e)]{Jarnicki})
$$\limsup_{\substack{z,w\to u\\z\neq
w}}\frac{b_D(z,w)}{|z-w|}\le\beta_D(u;1)$$ (cf. \cite[Lemma 4.3.3
(e)]{Jarnicki}) and the equality (see \cite[Remark, p. 11]{Jar-Nik})
$$\lim_{u\to\partial D}\beta_D(u;1)d_D(u)=\frac{\sqrt2}{2}.$$
\smallskip

Recall now another comparison result between $c_{D}$ and $k_{D}$
(see \cite[Proposition 9]{Nikolov2}): if $D$ is a finitely connected
bounded planar domain without isolated boundary points,\footnote{Any
$C^1$-smooth bounded planar domain is such a domain.} then
\begin{equation}\label{eq}
\lim_{\substack{w\to\partial D\\z\neq
w}}\frac{c_{D}(z,w)}{k_{D}(z,w)}=1
\quad\textup{uniformly in }z\in D.
\end{equation}
Similar results for $c_D,$ $k_D,$ $l_D$ and $b_D$ in the strictly
pseudoconvex case can be found in \cite[Theorem 1]{Venturini} and
\cite[Proposition 4]{Nikolov3}.

The next proposition shows that
\eqref{eq} remains true if we replace $c_D$ or $k_D$
by $b_D/\sqrt 2.$

\begin{proposition}
If $D$ is a finitely connected bounded planar domain without isolated boundary
points, then
$$\lim_{\substack{w\to\partial D\\z\neq w}}\frac{b_{D}(z,w)}{c_{D}(z,w)}=
\lim_{\substack{w\to\partial D\\z\neq
w}}\frac{b_{D}(z,w)}{k_{D}(z,w)}=\sqrt{2}
\quad\textup{uniformly in }z\in D.$$
\end{proposition}

\noindent{\bf Remark.} The isolated boundary points condition is
essential. Indeed, if $p$ is an isolated boundary point of a planar
domain $D\neq\Bbb C\setminus\{p\},$ then $c_D=c_{D\cup\{p\}}$ and
$b_D=b_{D\cup\{p\}},$ but $k_D(z,w)\to\infty$ as $w\to p$ and $z\in
D$ is fixed.

\section{Definitions}

\noindent{\bf 1.} A boundary point $p$ of a planar domain $D$ is said
to be Dini-smooth if $\partial D$ near $p$ is given by a Dini-smooth
curve $\gamma:[0,1]\rightarrow\mathbb{C}$ with $\gamma'\neq 0$ (i.e.,
$\ds\int_0^1\frac{\omega(t)}{t}dt<\infty,$ where $\omega$ is the
modulus of continuity of $\gamma'$). A planar domain is called
Dini-smooth if all its boundary points are Dini-smooth.
\smallskip

\noindent{\bf 2.} Let $D$ be a domain in $\mathbb{C}^{n}.$

The Bergman distance $b_D$ of $D$ is the integrated form of the Bergman
metric $\beta_D,$ i.e.,
$$b_D(z,w)=\inf_\gamma\int_0^1\beta_D(\gamma(t);\gamma'(t))dt,
\quad z,w\in D,$$
where the infimum is taken over all smooth curves $\gamma:[0,1]\to
D$ with $\gamma(0)=z$ and $\gamma(1)=w.$

Recall that
$$\beta_{D}(z;X)=\frac{M_{D}(z;X)}{K_{D}(z)},\quad z\in D,\
X\in\mathbb{C}^{n},$$ where
$$M_{D}(z;X) = \sup\{|f'(z)X| : f \in L^{2}_{h}(D),\
\lVert f\rVert_{L^2(D)}\le 1,\ f(z) = 0\}$$
and
$$K_{D}(z)=\sup\{|f(z)|:f\in L^{2}_{h}(D),\ \lVert f\rVert_{L^2(D)}\le 1\}$$
is the square root of the Bergman kernel on the diagonal (we assume
that $K_D>0;$ for example, this holds if $D$ is bounded).

The Carath\'eodory distance $c_{D}$ and the Lempert function
$l_{D}$ of $D$ are defined as follows:
$$c_{D}(z,w)=\sup\{\tanh^{-1}|f(w)|:f\in\mathcal{O}(D,\mathbb{D}),\,\hbox{ with }f(z)=0\},$$
$$l_{D}(z,w)=\inf\{\tanh^{-1}|\alpha|:\exists\varphi\in\mathcal{O}(\mathbb{D},D)
\hbox{ with }\varphi(0)=z,\varphi(\alpha)=w\},$$ where $\mathbb{D}$
is the unit disc.

The Kobayashi distance $k_{D}$ is the largest pseudodistance
not exceeding $l_{D}.$ It is well-known that $k_{D}$ is the
integrated form of Kobayashi metric $\kappa_{D}$ defined by
$$\kappa_{D}(z;X)=\inf\{|\alpha|:\exists\varphi\in\mathcal{O}(\mathbb{D},D)
\hbox{ with }\varphi(0)=z,\,\alpha\varphi'(0)=X\}.$$

Note that $c_D\le k_D\le l_D$ and $c_D\le b_D.$ On the other hand,
$k_D=l_D$ for any planar domain $D$ (cf. \cite[Remark
3.3.8(e)]{Jarnicki}).

We refer to \cite{Jarnicki} for other basic properties of the above
invariants.

\section{Proofs}

To prove Proposition 1, we shall need the following

\begin{lemma} (a)
$$\log\left(1+\frac{|z-w|}{\sqrt{(1-|z|^2)(1-|w|^2)}}
\right)\le\frac{b_{\Bbb D}(z,w)}{\sqrt 2}$$
$$\le\log\left(1+\frac{2|z-w|}{\sqrt{(1-|z|^2)(1-|w|^2)}}
\right);$$

\noindent (b) $$\log\left(1+\frac{|z-w|}{2\sqrt{d_{\Bbb D}(z)d_{\Bbb D}(w)}}\right)
\le\frac{b_{\Bbb D}(z,w)}{\sqrt 2}\le\log\left(1+\frac{\sqrt2|z-w|}
{\sqrt{d_{\Bbb D}(z)d_{\Bbb D}(w).}}\right).$$

\end{lemma}

\noindent{\it Proof.} (a) We have that
$$\sqrt2 b_{\Bbb D}(z,w)=2k_{\Bbb D}(z,w)=\log\frac{1+
\ds\left|\frac{z-w}{1-\bar zw}\right|}
{1-\ds\left|\frac{z-w}{1-\bar zw}\right|}=$$
$$\log\left(1+\frac{2|z-w|}{|1-\bar zw|-|z-w|}\right)=
\log\left(1+2|z-w|\frac{|1-\bar zw|+|z-w|}{(1-|z|^2)(1-|w|^2)}\right).$$
It remains to use that
\begin{equation}\label{for}
|1-\bar zw|^2=(1-|z|^2)(1-|w|^2)+|z-w|^2
\end{equation}
and hence
$$\sqrt{(1-|z|^2)(1-|w|^2)}\le|1-\bar zw|\le\sqrt{(1-|z|^2)(1-|w|^2)}+|z-w|.$$

\noindent (b) The lower estimate follows from (a) and $\ds d_{\Bbb D}(z)=1-|z|\ge\frac{1-|z|^2}{2}.$

To get the upper estimate, we have to show that
$$1+2|z-w|\frac{|1-\bar zw|+|z-w|}{(1-|z|^2)(1-|w|^2)}\le
\left(1+\frac{\sqrt2|z-w|}{\sqrt{(1-|z|)(1-|w|)}}\right)^2$$
which is equivalent to
$$\frac{|1-\bar zw|+|z-w|}{(1-|z|^2)(1-|w|^2)}\le
\frac{\sqrt2}{\sqrt{(1-|z|)(1-|w|)}}+\frac{|z-w|}{(1-|z|)(1-|w|)},$$
i.e.,
$$|1-\bar zw|\le\sqrt{2(1-|z|)(1-|w|)}(1+|z|)(1+|w|)+|z-w|(|z|+|w|+|zw|).$$
So, it is enough to prove that
$$|1-\bar zw|^2\le2(1-|z|)(1-|w|)(1+|z|)^2(1+|w|)^2+|z-w|^2(|z|+|w|+|zw|)^2.$$
Using \eqref{for} and dividing by $(1+|z|)(1+|w|),$
the last inequality
becomes
$$|z-w|^2(1-|z|-|w|-|zw|)\le(1-|z|-|w|+|zw|)(1+2|z|+2|w|+2|zw|)$$
which follows from $|z-w|^2\le |z|^2+|w|^2+2|zw|\le|z|+|w|+2|zw|.$
\smallskip

\noindent{\bf Remark.} (a) The constants 1 and 2 in front of $|z-w|$ in the lower
and upper estimates in Lemma 4 (a) are sharp. To see this, let
$$\frac{|z-w|^2}{(1-|z|^2)(1-|w|^2)}\to 0\mbox{ and }\infty,\mbox{ respectively}.$$

\noindent (b) The constants $\ds\frac{1}{2}$ and $\sqrt 2$ in front of $|z-w|$ in the lower
and upper estimates in Lemma 4 (b) are sharp, too. To see this, let $|z|\to 1$ and then 
$$\frac{|z-w|^2}{(1-|z|^2)(1-|w|^2)}\to 0\mbox{ and }w\to 0,\mbox{ respectively}.$$
\smallskip

\noindent{\it Proof of Proposition 1.} Let $D\supset
(z_n)\to p$ and $D\supset(w_n)\to q$ ($z_n\neq w_n$). It is enough to find a
constant $c>1$ such that the respective estimates for $b_D(z_n,w_n)$
hold for any $n.$

Note that, by \cite[Proposition 5 and Corollary 6]{Nikolov2}, for
any neighborhood $U$ of $p$ there exist a neighborhood $V\subset U$
and a constant $c_1>0$ such that
\begin{equation}\label{main}
|{\sqrt2}b_D(z,w)+\log d_D(z)+\log d_D(w)|<c_1,\quad z\in D\cap V, w\in D\setminus U.
\end{equation}

This inequality provides the desired constant if $D\ni p\neq q\in
D,$ or $p\in\partial D,$ $q\in D,$ or $p\in D,$ $q\in\partial D,$ or
$\partial D\ni p\neq q\in\partial D.$

For a planar domain $\Omega,$ set
$\beta_\Omega(z)=\beta_\Omega(z;1),$ $M_\Omega(z)=M_\Omega(z;1)$ and
$\kappa_\Omega(z)=\kappa_\Omega(z;1).$

If $p=q\in D,$ then the continuity of $\beta_D$ implies that
$$\frac{b_D(z_n,w_n)}{|z_n-w_n|}\to\beta_D(p)>0$$
and we may easily find the desired constant.

It remains to consider the most difficult case $p=q\in\partial D.$
Some of our arguments will be close to that in the proof of
\cite[Proposition 5]{Nikolov2}.

This proof allows us to assume that $p=1$ and
$$\{z\in\Bbb D:|z-1|<r\}=:E_r\subset D\subset\Bbb D$$
for some $r>0$ (after an appropriate conformal map). Then
\begin{equation}\label{ber}\sqrt2\frac{\kappa^2_{\Bbb D}(z)}
{\kappa_{E_r}(z)}=\frac{M_{\Bbb D}(z)}{K_{E_r}(z)}\le\beta_D(z)\leq
\frac{M_{E_r}(z)}{K_{\Bbb D}(z)}=
\sqrt2\frac{\kappa^2_{E_r}(z)}{\kappa_{\Bbb D}(z)},\quad z\in E_r
\end{equation}
(the both equalities hold because $E_r$ is a simply connected
domain).

Fix an $r_1\in (0,r).$ The localization of the Kobayashi metrics
from \cite[Theorem 2.1 and Lemma 2.2]{Forstneric} implies that
\begin{equation}\label{kob}
\kappa_{\Bbb D}(z)>(1-c_2d_{\Bbb D}(z)) \kappa_{E_r}(z),\quad z\in
E_{r_1},
\end{equation} for some constant $c_2>0.$
Choose an $r_2\in(0,r_1]$ with $2c_2r_2\le 1.$ Then
$$(1-c_2d_{\Bbb D}(z))\kappa_{\Bbb D}(z)<\frac{\beta_D(z)}{\sqrt 2}
<(1+2c_2d_{\Bbb D}(z))\kappa_{\Bbb D}(z),\quad z\in
E_{r_2}.$$ Since $\ds\kappa_{\Bbb D}(z)=\frac{\beta_{\Bbb
D}(z)}{\sqrt 2}=\frac{1}{1-|z|^2},$ it follows for $\ds
c_3=2\sqrt 2c_2$ that
\begin{equation}\label{est}
\frac{1}{\sqrt2}\left(\frac{1}{d_{\Bbb D}(z)}-c_2\right)<\beta_D(z)<\beta_{\Bbb D}(z)+c_3,\quad z\in E_{r_2}.
\end{equation}

We may assume that $z_n,w_n\in E_{r_3},$ where $r_3\in(0,r_2/2]$ is
such that if $\alpha_n$ is the shorter arc with endpoints $z_n$ and
$w_n$ of the circle through $z_n$ and $w_n$ which is orthogonal to
the unit circle, then $\alpha_n\subset E_{r_2}.$ Hence
$$b_{D}(z_n,w_n)<\int_{\alpha_n}\left(\frac{\sqrt 2}{1-|z|^2}+c_3\right)dl$$
$$=b_{\Bbb D}(z_n,w_n)+c_3l(\alpha_n)<b_{\Bbb D}(z_n,w_n)+2c_3|z_n-w_n|$$
for any $n.$

Now, using Lemma 4 (b) and the equality
\begin{equation}\label{ine}
d_{\Bbb D}(z)=d_D(z),\quad z\in E_{r_3},
\end{equation}
it is easy to find a constant $c>1$ such that the upper estimate for
$b_D(z_n,w_n)$ in Proposition 1 holds for any $n.$

It is left to manage the lower estimate. Let $\gamma_n:[0,1]\to D$
be a smooth curve such that $\gamma_n(0)=z_n,$ $\gamma_n(1)=w_n$ and
$$b_D(z_n,w_n)+|z_n-w_n|>\int_0^1\beta_D(\gamma_n(t);\gamma_n'(t))dt.$$

Consider the set $A$ of all $n$ for which $\gamma_n(0,1)\not\subset E_{r_2}.$
For any $n\in A$ we may find a number $t_n\in(0,1)$ such that
$|u_n-1|=r_2,$ where $u_n=\gamma(t_n).$ By \eqref{main}, there
exists a constant $c_4>0,$ which does not depend on $n\in A,$ such
that
$$b_D(z_n,w_n)+|z_n-w_n|>b_D(z_n,u_n)+b_D(u_n,w_n)$$
$$>-\frac{\log d_D(z_n)}{\sqrt 2}-\frac{\log d_D(w_n)}{\sqrt
2}-c_4.$$ This inequality easily provides a constant $c>1$ for
which the lower estimate for $b_D(z_n,w_n)$ in Proposition 1
holds for any $n\in A.$

Let now $n\not\in A.$
Since
$$d_{\Bbb D}(\gamma_n(t))\le f_n(t):=d_{\Bbb D}(z_n)+|z_n-\gamma_n(t)|<2r_3+r_2\le 2r_2\le 1/c_2,$$
$$d_{\Bbb D}(\gamma_n(t))\le g_n(t):=d_{\Bbb D}(w_n)+|w_n-\gamma_n(t)|<1/c_2$$ and $|s|'\ge|s'|,$
it follows by \eqref{est} that, for any $t_n\in(0,1),$
$$\sqrt{2}(b_D(z_n,w_n)+|z_n-w_n|)>\int_0^1\left(\frac{1}{d_{\Bbb D}(\gamma_n(t))}-c_2\right)|\gamma_n'(t)|dt$$
$$\ge\int_0^{t_n}\left(\frac{1}{f_n(t)}-c_2\right)df_n(t)-\int_{t_n}^1\left(\frac{1}{g_n(t)}-c_2\right)dg_n(t)$$
$$=\log\left(1+\frac{|z_n-\gamma_n(t_n)|}{d_{\Bbb D}(z_n)}\right)-c_2|z_n-\gamma_n(t_n)|$$
$$+\log\left(1+\frac{|w_n-\gamma_n(t_n)|}{d_{\Bbb D}(w_n)}\right)-c_2|w_n-\gamma_n(t_n)|$$
$$>\log\left(1+\frac{|z_n-\gamma_n(t_n)|.|w_n-\gamma_n(t_n)|}{c_5d_{\Bbb D}(z_n)d_{\Bbb D}(w_n)}\right)$$
for some constant $c_5>1.$ Choosing now $t_n$ such that $|z_n-\gamma_n(t_n)|=|w_n-\gamma_n(t_n)|$
and using \eqref{ine}, we obtain the lower estimate in Proposition 1.

So, Proposition 1 is completely proved.
\smallskip

\noindent{\it Proof of Proposition 3.} By the K\"obe uniformization
theorem, we may assume that $\partial D$ consists of disjoint
circles. Using Proposition 1, Corollary 2, \eqref{eq} and
compactness, it is enough to prove that
$$\lim_{\substack{z,w\to p\\z\neq w}}\frac{b_D(z,w)}{k_D(z,w)}=\sqrt
2$$ for any point $p\in\partial D.$

Applying an inversion, we may assume that the outer boundary of $D$
is the unit circle and $p=1.$ Then \eqref{ber} and
\eqref{kob} imply
$$\lim_{z\to 1}\frac{\beta_{E_r}(z)}{\beta_D(z)}=1=
\lim_{z\to 1}\frac{\kappa_{E_r}(z)}{\kappa_D(z)}.$$

The first equality shows that $\ds\liminf_{\substack{z,w\to
1\\z\neq w}}\frac{b_{E_r}(z,w)}{b_D(z,w)}\ge 1.$

To get that
\begin{equation}\label{inf}
\limsup_{\substack{z,w\to 1\\z\neq
w}}\frac{b_{E_r}(z,w)}{b_D(z,w)}\le 1,
\end{equation}
we shall follow the proof of \cite[Proposition 3]{Venturini}. Fix an
$\eps>0$ and choose an $r_1\in(0,r)$ such that
$$\beta_{E_r}(z)<(1+\eps)\beta_D(z),\quad z\in E_{r_1}.$$
Combining the argument in the case $n\not\in A$ from the previous
proof and the estimates from Proposition 1, we may find an
$r_2\in(0,r_1)$ such that if $z,w\in E_{r_2}$ and $\gamma:[0,1]\to
D$ is a smooth curve for which $\gamma(0)=1,$ $\gamma(1)=w$ and
$$\int_0^1\beta_D(\gamma(t);\gamma'(t))dt\le(1+\eps)b_D(z,w),$$
then $\gamma([0,1])\subset E_{r_1}.$ It follows that
$$b_{E_r}(z,w)\le\int_0^1\beta_{E_r}(\gamma(t);\gamma'(t))dt$$
$$\le (1+\eps)\int_0^1\beta_D(\gamma(t);\gamma'(t))dt\le
(1+\eps)^2b_D(z,w),\quad z,w\in E_{r_2}.$$ To obtain \eqref{inf}, it
remains to let $\eps\to 0.$

So, $\ds\lim_{\substack{z,w\to 1\\z\neq
w}}\frac{b_{E_r}(z,w)}{b_D(z,w)}=1.$

On the other hand, \cite[Proposition 8]{Nikolov2} gives the
estimates from Proposition 1 for $2k_D$ instead of $\sqrt 2 b_D.$
Then we obtain as above
$$\lim_{\substack{z,w\to 1\\z\neq
w}}\frac{\kappa_{E_r}(z,w)}{\kappa_D(z,w)}=1.$$

Now, the equality $b_{E_r}=\sqrt2k_{E_r}$ completes the proof.


\begin{thebibliography}{}

\bibitem{Bonk} Z. M. Balogh, M. Bonk, {\it Gromov hyperbolicity
and the Kobayashi metric on strictly pseudoconvex domains},
Comment. Math. Helv. 75 (2000), 504--533.

\bibitem{Forstneric} F. Forstneric, J.-P. Rosay
{\it Localization ot the Kobayashi metric and the boundary continuity
of proper holomorphic mappings}, Math. Ann. 279 (1987), 239--252.

\bibitem{Herbort} G. Herbort, {\it Estimation of the
Carath\'eodory distance on pseudoconvex domains of finite type,
whose boundary has a Levi form of corank at most one},
Ann. Polon. Math. 109 (2013), 209--260.

\bibitem{Jar-Nik} M. Jarnicki, N. Nikolov,
{\it Behavior of the Carath\'eodory metric near strictly convex
boundary points}, Univ. Iag. Acta Math. XL (2002), 7--12.

\bibitem{Jarnicki} M. Jarnicki, P. Pflug, {\it Invariant distances and
metrics in complex analysis}, de Gruyter Exp. Math. 9, de Gruyter, Berlin.

\bibitem{Nikolov2} N. Nikolov, {\it Estimates of invariant distances
on "convex" domains}, Ann. Mat. Pura Appl. 193 (2014), 1595--1605.

\bibitem{Nikolov3} N. Nikolov, {\it Comparison of invariant functions on
strongly pseudoconvex domains}, J. Math. Anal. Appl. 421 (2015), 180--185.

\bibitem{Nikolov1} N. Nikolov, P. Pflug, P. J. Thomas, {\it Upper bound
for the Lempert function of smooth domains}, Math. Z. 266 (2010),
425--430.

\bibitem{Venturini} S. Venturini, {\it Comparison between the Kobayashi
and the Carath\'{e}odory distances on strongly pseudoconvex bounded
domains in $\mathbb{C}^{n}$}, Proc. Am. Math. Soc. 107 (1989),
725--730.

\end{thebibliography}
\end{document}